\documentclass[12pt,leqno]{article}

\usepackage{amsmath}
\usepackage{amssymb}

\def\a{\alpha}

\def\p{\partial}

\def\s{\sigma}

\begin{document}

\title{Some remarks on Myers theorem for Finsler manifolds}
\author{Mihai Anastasiei}
\date{}

\maketitle

\begin{abstract}

The standard Bonnet-Myers theorem says that if the Ricci scalar of a Riemannian manifold is bounded  below by a positive number, then the manifold is compact. Moreover, a bound of its diameter is pointed out. The theorem was extended to Finsler manifolds. In this paper we prove  that if a certain condition on the average of the Ricci scalar holds, then the Finsler manifold $M$ is compact if the Ricci scalar is bounded above by the same positive number. An upper bound of the diameter is also found. With no condition on Ricci scalar itself but with a different one on its average, we find that the Finsler manifold $M$ is again compact. This time no bound of the diameter is found.The proofs are given in the Finslerian setting and are based on the index form along geodesics.


\end{abstract}

\section*{Introduction}

In the textbook \cite{BCS} by D. Bao, S.-S.Chern, Z. Shen a detailed  proof of the following version of the classical Myers theorem is provided.

{\bf Bonnet-Myers Theorem .} {\it Let $(M,F)$ be a forward geodesically complete connected
Finsler manifold of dimension $n$. Suppose its Ricci scalar $Ric$ has
the following uniform positive lower bound
$$
Ric\geq (n-1)a>0.
$$

Then:

\begin{itemize}
\item[$(1)$] Along every geodesic the distance between any two
successive conjugate points is at most $\dfrac{\pi}{\sqrt{a}}$.

\item[$(2)$] The diameter of $M$ is at most $\dfrac{\pi}{\sqrt{a}}$.

\item[$(3)$] $M$ is in fact compact.

\item[$(4)$] The fundamental group $\pi(M,x)$ is
finite.

\end{itemize}}

Carefully reading the said proof  in \cite{BCS}, p. 194-198
it comes out that the essential step is to prove its first statement.

We reformulate this first statement  as follows.

{\bf Theorem 1} {\it Let $\s(t)$, $0\leq t \leq L$ be a unit speed
geodesic with velocity field $T$ and $Ric(t):= Ric_{(\s(t),T)}.$ If

\begin{itemize}
\item[a)] The Ricci scalar $Ric(t)$ satisfies
$$Ric(t) \geq (n-1)a>0 ,$$ for a constant $a>0$ and for every $t\in [0,L]$ and

\item[b)] $$L\geq \frac{\pi}{\sqrt{a}},$$

\end{itemize}

then $\s$ must contain  conjugate points to $\s(0)$.}

In \cite{AN1} we have presented, again in the Finslerian setting, a generalization of the Bonnet-Myers theorem allowing and negative values for $Ric$ . This fact is important for General Relativity, \cite{Ga1}. Again the first statement is essential. We restate it as follows.

{\bf Theorem 2.} {\it Let $\s(t)$, $0\leq t \leq L$ be a unit speed
geodesic with velocity field $T$ and $Ric(t):= Ric_{(\s(t),T)}.$ If

\begin{itemize}
\item[a)] The Ricci scalar $Ric(t)$ satisfies
$$Ric(t) \geq (n-1)a + \frac{df}{dt} ,$$ for some function $f$ with $|f(t)|\leq \frac{\Lambda}{\pi}, \Lambda \geq 0,$ and
\item[b)]
$$ L\geq \frac{\Lambda}{a(n-1)}+ \sqrt{\frac{\pi^2}{a}+\frac{\Lambda ^2}{a^2(n-1)^2}},$$

\end{itemize}

then $\s$ must contain  conjugate points to $\s(0)$.}

If the Finsler manifold  $(M,F)$ is forward geodesically complete and connected, by the Theorem 2 it follows that it is compact and its diameter satisfies
$$ diam(M) \leq \frac{\Lambda}{a(n-1)}+ \sqrt{\frac{\pi^2}{a}+\frac{\Lambda ^2}{a^2(n-1)^2}}.$$

For $\Lambda =0$, the function $f$ must be zero and the Theorem 2 reduces to the Theorem 1. It is also clear that the Theorem 2 brings a novelty only when along of the geodesic we have
$$(n-1)a+ \frac{df}{dt} \leq Ric(t) < (n-1)a,\leqno(*) $$ that is when the function $f$ is decreasing. Thus negative values for $Ric$ are allowed. For instance, if we take $f(t) =\frac{\Lambda}{\pi}\cos\pi t,$ for $t=\frac{1}{2}$ and $\Lambda > (n-1)a$ the left side of (*) is negative.

In this paper we prove

{\bf Theorem A.} {\it Let $(M,F)$ be a forward geodesically complete connected
Finsler manifold of dimension $n$. Suppose that
\begin{itemize}
\item[a)] The Ricci scalar $Ric$ has the following uniform positive upper bound
$$
Ric < (n-1)a
$$
for a constant $a>0,$
\item[b)] For every geodesic $\sigma$ parameterized by the arc-length $t\in [0,L]$ we have $$\int_{0}^{L}Ric(t)dt\geq a(n-1)L + \varepsilon\Lambda ,
 $$
 for $\varepsilon = \pm 1$ and a constant $\Lambda > 0.$
\end{itemize}

Then:

\begin{itemize}
\item[$(1)$] Along every geodesic the distance between any two
successive conjugate points is at most $-\varepsilon\frac{\Lambda}{a(n-1)}+ \sqrt{\frac{\pi^2}{a}+\frac{\Lambda ^2}{a^2(n-1)^2}}$.

\item[$(2)$] The diameter of $M$ is at most $-\varepsilon\frac{\Lambda}{a(n-1)}+ \sqrt{\frac{\pi^2}{a}+\frac{\Lambda ^2}{a^2(n-1)^2}}$.

\item[$(3)$] $M$ is  compact.

\end{itemize}}

In the assumption b) from the hypothesis of the Theorem A the arc-length parameter runs from $0$ to $L$. But if the manifold is forward geodesically complete every  such geodesic can be extended to a geodesic defined on $[0,\infty)$. On such an extended  geodesic the assumption b)  should be replaced with the assumption that the integral of $Ric$ on $[0,\infty)$ diverges to $+\infty$. What happens in this case ? An answer is given by

{\bf Theorem B.} {\it Let $(M,F)$ be a forward geodesically complete connected Finsler  manifold of dimension $n$. If there exists a point $p\in M$ such that along each geodesic $\sigma:[0,\infty)\to M$ emanating from $p$ and parameterized by arc length $t$ the condition

$$
\int^\infty_0  Ric(t)dt=\infty,
$$
holds, then $M$ is compact.}

Now the hypothesis say nothing about $Ric$ itself but we are not able to evaluate the diameter of the manifold.

All the results above mentioned and the others of this type essentially depend on the existence of a pair of conjugate points along the geodesics of a Riemannian, Finslerian or Lorentzian manifold. The problem of the said existence has been appeared also in General relativity either in order to prove a singularity theorem (Hawking and Penrose,\cite{HP}) or to argue a relationship between energy density and the size of space( Galloway, \cite{Ga1}, Frankel and Galloway,\cite{FG}). Various aspects of this problem remained of interest along years (Chicone and Ehrlich,\cite{CE}. Kim and Kim,\cite{KK}). Many of them have been considered in more general spacetimes (Caponio,Javaloyes and Masillo,\cite{CJM}, Galloway and Woolgar,\cite{GW}) and some of them  are also of interest for the geometric structures of the extended theories of gravity studied by S. Capozziello and his co-workers ( \cite{Ca1},\cite{Ca2})  or by S. Vacaru (\cite{Va1},\cite{Va2}).

The structure of the paper is as follows. In the first Section we recall the results from Finsler geometry to be used. The textbook \cite{BCS} will be closely followed. The next two Sections are devoted to the proofs of the Theorem A and the Theorem B, respectively.

\section{Preliminaries}

We shall use the notations,
the terminology and results from \cite{BCS} without comments.

Let $(M,F)$ be a Finsler manifold. The Finsler structure $F$ is a
function $F:TM\to [0,\infty)$, $(x,y)\to F(x,y)$ which is
$C^{\infty}$ on the slit tangent bundle $TM\backslash 0$,
positively homogeneous in $y$ and whose Hessian matrix $g_{ij}:=
\dfrac{1}{2}\dfrac{\p^{2}F^{2}}{\p y^{i}\p y^{i}}$ is
positive-definite at every point of $TM\backslash 0.$

The Chern connection of local coefficients $\Gamma ^i_{jk}(x,y)$ is a linear connection in the pull-back
bundle $\pi^{*}TM$ over $TM\backslash 0$, where $\pi:TM\to M$ is
the natural projection. It is only $h$-metrical and it has two
curvatures $R_{j} \ ^{i}\ _{kh}$, $P_{j} \ ^{i}\  _{kh}$.

Let be $y$ a non zero element of $T_{x}M$. Then $g(x,y) =
g_{ij}(x,y)dx^{i}\otimes dx^{j}$ is an inner product which is used
to measure lengths and angles in $T_{x}M$.One calls $y$ a flagpole
of the flag (a plane in $T_{x}M$) spanned by $l=\dfrac{y}{F(x,y)}$,
and another unit vector $V$ which is orthogonal to the flagpole.

The flag curvature is then given as
$$
K(x,y,l\wedge V):= V^{i}(l^{j}R_{jikh}l^{h})V^{k}=:
V^{i}R_{ik}V^{k}. \leqno(1.1)
$$
 Let $\{l, e_{\a},
\a=1,\ldots,n-1\}$ be a $g$-orthonormal basis for the fiber of
$\pi^{*}TM$ over the point $(x,y)\in TM\backslash 0.$ With respect
to it one has $K(x,y,l\wedge e_{\a})=R_{\a\a}$. The Ricci scalar
denoted by Ric is
$$
Ric:= \sum_{\a=1}^{n-1}K(x,y,l\wedge e_{\a}) =
\sum_{\a=1}^{n-1}R_{\a\a}. \leqno(1.2)
$$

If $(M,F)$ has constant flag curvature $c$, then
$$
Ric = (n-1)c. \leqno(1.3)
$$

Let $\s (t), 0\leq t \leq L,$  be a unit geodesic with velocity
field $T$. One abbreviates $g_{(\s,T)}$ by $g_{T}$.

For a vector field $W(t):= W^{i}(t)\dfrac{\p}{\p x^{i}}$ along
$\s$, the expression,
$$
D_{T}W =
\left[\dfrac{dW^{i}}{dt}+W^{j}T^{k}(\Gamma_{jk}^{i}(\s,T))\right]\dfrac{\p}{\p
x^{i}}\leqno(1.4)
$$
is called covariant derivative with reference vector $T$.

The constant speed geodesics are solutions of $D_{T}T=0,$ with
reference vector $T$.

One says that $W$ is parallel long $\s$ if $D_{T}W=0,$ with
reference vector $T$. Parallel transport (with reference vector
$T$) one defines on the standard way. The parallel
transport preserves $g_{T}$-lengths and angles.

For two continuous and piecewise $C^{\infty}$ vector fields $V$
and $W$ along $\s$ the index form is
$$
I(V,W) = \int_{0}^{L}[g_{T}(D_{T}V, D_{T} W) -
g_{T}(R(V,T)T,W)]dt. \leqno(1.5)
$$
It can be re-expressed in the form
$$
I(V,W) = \int_{0}^{L}[g_{T}(D_{T}V, D_{T} W) -
K(T,W)g_{T}(W,W)]dt, \leqno(1.5')
$$
where $K(T,W)$ is the flag curvature of the flag with flagpole $T$ and transverse edge $W$.

Here all $D_{T}$ are calculated with reference vector $T$ and
$$R(V,T)T:=(T^{j}R_{jkh}^{i}T^{h})V^{k}\frac{\p}{\p x^{i}}$$ is
evaluated at the point $(\s,T)$.

The index form is bilinear and symmetric.

Let $0=:t_0<t_1<...<t_h:=L$ be a partition of $[0,L]$ such that $V$ and $W$ are both $C^\infty$ on each closed subinterval $[t_{s-1},t_s]$. Using integration by parts, one can rewrite the index form as
\begin{equation}
\begin{array}{ll}
I(V,W):=& g_T(D_T V,W)\Bigg{|}^{L}_{0}-\left.\displaystyle\sum^{h-1}_{s=1}g_T(D_T V,W)\right|^{t_s^+}_{t^-_s}-\\ \\
& -\displaystyle\int^L_0 g_T(D_TD_TV+R(V,T)T,W)dt.
\end{array}
\tag{1.6}
\end{equation}

A vector field $J$ along $\sigma$ is said to be a {\it Jacobi field} if it satisfies the equation
$$
D_TD_TJ+R(J,T)T=0.\leqno(1.7)
$$

One says that $q=\sigma(L)$ is conjugate with $p=\sigma(0)$ along $\sigma$ if there exists a nonzero Jacobi field $J$ along $\sigma$ which vanishes at $p$ and $q$ i.e. $J(0)=J(L)=0$.
We recall  from \cite{BCS} p.182 the following

{\bf Proposition 1.1}  {\it Let $\s(t), 0 \leq t\leq
r$ be a geodesic in a Finsler manifold $(M,F)$. Suppose no point
$\s(t), 0 < t\leq r$ is conjugate to $p:=\s(0)$. Let $W$ be any
piecewise $C^{\infty}$ vector field along $\s$ and let $J$ denote
the unique Jacobi field along $\s$ that has the same boundary
values as $W$. That is, $J(0) = W(0)$ and $J(r) = W(r)$. Then}
$$
I(W,W)>I(J,J). \leqno(1.8)
$$
{\it Equality holds if and only if $W$ is actually a Jacobi field, in
which case the said $J$ coincides with $W$.}

\section{Proof of Theorem A}

It suffices to prove that if along every unit speed
geodesic $\s(t)$, $0\leq t \leq L$ the Ricci scalar satisfies the hypothesis a) and b) of the Theorem A and if
$$ L\geq - \varepsilon\frac{\Lambda}{a(n-1)}+ \sqrt{\frac{\pi^2}{a}+\frac{\Lambda ^2}{a^2(n-1)^2}},$$
then $\s$ must contain  conjugate points to $\s(0).$

Using the parallel transport with reference vector
$T$ we construct a moving frame $\{e_{i}(t)\}$ along $\s$ such
that

(i) Each $e_{i}$ is parallel along $\s$, that is $D_{T}e_{i} = 0,$

(ii) $\{e_{i}(t)\}$ is a $g_{T}$-orthonormal frame,

(iii) $e_{n}= T.$

Define $W_{\a}(t) = f(t)e_{\a}(t)$ for some smooth function
$f$, $\a=1,2$, ..., $n-1$.

Fix a positive $r\geq L$ and consider the index from $I$ for
$\s(t), 0 \leq t \leq r$. By (1.5') we have
$$
I(W_{\a}, W_{\a}) = \int_{0}^{r}[g(D_{T}W_{\a},D_{T}W_{\a}) -
g(W_{\a},W_{\a})K(T,W_{\a})] dt,
$$
where
$K(T,W_{\a})$ is the flag curvature evaluated at the point
$(\s(t),T)\in TM\backslash 0.$

We have $D_{T}W_{\a}=\frac{df_{\a}}{dt}e_{\a}$ and
since the flag
curvature does not depend on vectors spanning  the flag, the equality
$K(T,W_{\a}) = K(T,e_{\a})$ holds.

Using these facts, $I(W_{\a},W_{\a})$ takes the form
$$
I(W_{\a},W_{\a}) =
\int_{0}^{r}\left[\left(\frac{df}{dt}\right)^{2} -
f^{2}K(T,e_{\a})\right]dt.\leqno(2.1)
$$
We take $f(t) = \sin\dfrac{\pi t}{r}$ and we get
$$
I(W_{\a},W_{\a}) = \frac{\pi^{2}}{2r} -
\int_{0}^{r}\sin^{2}\frac{\pi t}{r}K(T,e_{\a})dt.\leqno(2.2)
$$

Summing over $\a$ one obtains
$$
\sum_{\a}I(W_{\a},W_{\a}) = (n-1)\frac{\pi^{2}}{2r} -
\int_{0}^{r}Ric (t)dt+\int_{0}^{r}Ric (t)\cos^{2}\frac{\pi t}{r}dt.\leqno(2.3)
$$
By the assumptions a) and b) one gets
$$
\sum_{\a}I(W_{\a},W_{\a})\leq (n-1)\frac{\pi^{2}}{2r} -
(n-1)ar -\varepsilon \Lambda + (n-1)a\int_{0}^{r}cos^{2}\frac{\pi t}{r}dt\leqno(2.4)
$$
Computing the indicated integral one yields
$$
\sum_{\a}I(W_{\a}, W_{\a}) \leq \frac{(n-1)}{2r}(\pi ^2 - 2\varepsilon\frac{\Lambda}{n-1}r - ar^2)\leqno(2.5)
$$

and we have $\sum_{\a}I(W_{\a},W_{\a})\leq 0$ if $r\geq L=-\varepsilon\frac{\Lambda}{a(n-1)}+ \sqrt{\frac{\pi^2}{a}+\frac{\Lambda ^2}{a^2(n-1)^2}}.$ It follows that some $I(W_{\a},W_{\a})$ must
be non-positive and let denote that $W_{\a}$ by $W$.

We proceed by contradiction. Suppose that $\s(t), 0 \leq t\leq
r=- \varepsilon\frac{\Lambda}{a(n-1)}+ \sqrt{\frac{\pi^2}{a}+\frac{\Lambda ^2}{a^2(n-1)^2}}$ contains no conjugate
points. By the very definition of the conjugate points, the unique Jacobi field which vanishes at
the endpoints of $\s(t), 0\leq t\leq r$ is identically zero.
The vector field $W$ satisfies  $W(0) = W(r) = 0$
and it can not be a Jacobi field since is nowhere zero on $(0,r)$. By the Proposition 1.1 we have $0=I(J,J) <I(W,W)\leq 0$
which is a contradiction.

In order to prove  the statements 2)-3) of the Theorem A  the same arguments as those from \cite{BCS} p. 196-198 are used. We outline them in the following.

Since $M$ is forward geodesically complete, by the
Hopf-Rinow theorem any pair of points in $M$ can be joined by a
minimal geodesic.It is known that the cut point of $\s(0)$ appears before or coincide with the first conjugate point to $\s(0)$. As we have just proved, such a geodesic must have the length less than or equal with $-\varepsilon \frac{\Lambda}{a(n-1)}+ \sqrt{\frac{\pi^2}{a}+\frac{\Lambda ^2}{a^2(n-1)^2}}$ .
Thus $\mbox{diam}\, (M)\leq -\varepsilon\frac{\Lambda}{a(n-1)}+ \sqrt{\frac{\pi^2}{a}+\frac{\Lambda ^2}{a^2(n-1)^2}}$ , hence 2) holds.
By the statement 2) the manifold $M$ is forwardly bounded from the above. As it is always closed in its own topology, using again the Hopf-Rinow theorem one
concludes that $M$ is compact, that is, the statement 3) holds.Thus the Theorem A is completely proved. $\hfill\square$

{\bf Remark 2.1.} If in the main theorem (Theorem 1.2) from \cite{Wu} the function $max$ is explicitly written, two statements are obtained. The one covers the first three items of the Bonnet-Myers theorem. The other one is similar with the case $\varepsilon=-1$ from  Theorem A except that the bound of the diameter of $M$ is $\frac{\pi}{\sqrt{a}} + \frac{\Lambda}{a(n-1)}$. This is clearly lesser then our bound in the case $\varepsilon=-1$. But our bound in the case $\varepsilon=+1$  is strictly lesser then the bound $\frac{\pi}{\sqrt{a}} + \frac{\Lambda}{a(n-1)}$. The latter was  found in \cite{Wu} by using a Ricatti inequation satisfied by the trace of the Hessian of the Finslerian distance function on $M$. Thus we have three different bounds for the diameter of $M$  all depending on  $\Lambda$. If $\Lambda $ increase to $+\infty$ two of them monotonically increase also to $+\infty$ and one monotonically decreases to zero. For $\Lambda =0$ all three reduce to the bound given by the Bonnet-Myers theorem.

\section{Proof of Theorem B}

Before going on we notice that in the proof of the Theorem A a main fact was that for given a point $p\in M$ every unit speed geodesic emanating from $p$ contains a first point conjugates to $p$. Then using the Morse index form a evaluation of length of the geodesic from $p$ to this first conjugate point was performed. Based on it a bound of the diameter of $M$ was found and from here the conclusion that $M$ is compact. But the same conclusion can be derived directly from the just mentioned main fact. In the Riemannian case the remark is due to W. Ambrose (\cite{Am}). In our framework it can be formulated as follows.

{\bf Lemma 3.1.} {\it Let $(M,F)$ be a forward geodesically complete connected Finsler  manifold of dimension $n$. If there exists a point $p\in M$ such that every geodesic ray emanating from $p$ has a point conjugate to $p$ along that ray, then $M$ is compact.}

{\bf Proof.} Let $S_p$ be the indicatrix in the point $p\in M$. For each $y\in S_p$ issue the unit speed geodesic from $p$ with the initial unity velocity $y$. Each such geodesic is defined for any $t \in [0,\infty).$ Let $c_y$ be the value of $t$ in the first conjugate point of $p$ and $i_y$ the value of $t$ in the cut point of $p$. By the hypothesis of the Lemma 3.1 the set of $c_y$ is forwardly bounded from above (if $c_y =\infty$ one says that $p$ has no conjugate points along that geodesic) and since one has $i_y\leq c_y$ it follows that $ \sup _{y\in S_p} i_y\leq \sup _{y\in S_p}c _y$ and because the diameter of $M$ is less or equal with $ \sup _{y\in S_p} i_y $ it comes out that $M$ is forwardly bounded from the above. As $M$ is closed in its own topology, by the Hopf-Rinow theorem it is compact. $\hfill\square$

Thus in order to prove the Theorem B it suffices to  prove that there exists a point $p\in $ such that every unit speed geodesic $\sigma:[0,\infty)\to M$ issuing from $p$ has a point conjugate to $p$ along $\sigma:[0,\infty)\to M$ . The Morse index lemma will be again used.
We repeat the construction leading to the formula (2.1) from Section 2 and replace the function $f$ by the following one: $f(t)= 1 , t\in [0,1], = 1, t\in [1,b] , = \frac{r-t}{r-b} , t\in [b,r]. $ for any $b\in (1,r)$. Then summing over $\a$, instead of (2.3) one gets

$$
\sum_{\a}I(W_{\a},W_{\a}) = \int_0^1 ((n-1)- t^2Ric(t))dt - \int_1^b Ric(t)dt+$$

$$
+ \int_b^r ( \frac{n-1}{(r-b)^2}- \frac{(r-t)^2}{(r-b)^2}Ric(t))dt .
$$

In the right hand of this equality, the first integral is finite, by the hypothesis of the Theorem B the second integral in (3.1) diverges to $-\infty$ and by an integration by parts it comes out that the third integral tends to $0$ when $r$ tends to $\infty$.

Thus $\sum_{\a}I(W_{\a},W_{\a})\leq 0$ and as we have seen before this fact implies that $p$ has a conjugate point along the geodesic $\sigma:[0,\infty)\to M$, q.e.d.

{\bf Acknowledgments} The author was partially supported by a Grant of the Romanian National
Authority for Scientific Research, CNSS-UEFISCDI, project number PN-II-IDPCE-
2011-3-0256.

{\it Author's address:}\\
Mihai Anastasiei\\
Faculty of Mathematics,\\
Alexandru Ioan Cuza University of Ia\c si\\
and\\
Mathematical Institute "O. Mayer", Romanian Academy, \\
Ia\c{s}i, Romania,\\
E-mail: anastas@uaic.ro.

\end{document}